\documentclass{amsart}

\usepackage{amssymb} \usepackage{amsfonts} \usepackage{amsmath}
\usepackage{amsthm} \usepackage{epsfig} 
\usepackage{color}
\usepackage{amscd}
\usepackage[all]{xy}

\newcommand{\dimo}[1]{\vspace{2pt}\noindent\textit{Proof of {#1}}.\ }
\newcommand{\finedimo}{{\hfill\hbox{$\square$}\vspace{2pt}}}

\newtheorem{lemma}{Lemma}
\newtheorem{teo}[lemma]{Theorem}
\newtheorem{prop}[lemma]{Proposition}
\newtheorem{cor}[lemma]{Corollary}

\theoremstyle{definition}

\newtheorem{oss}[lemma]{Remark}
\newtheorem{quest}[lemma]{Question}

\theoremstyle{remark}

\newcommand{\Vol}{{\rm Vol}}

\newcommand{\calS}{\ensuremath {\mathcal{S}}}

\newcommand{\matZ} {\ensuremath {\mathbb{Z}}}

\author{Bruno Martelli}
\address{Dipartimento di Matematica ``Tonelli'', Largo Pontecorvo 5, 56127 Pisa, Italy}
\email{martelli at dm dot unipi dot it}

\title{Hyperbolic three-manifolds \\ that embed geodesically}

\begin{document}

\begin{abstract}
We prove that every complete finite-volume hyperbolic 3-manifold $M$ that may be tessellated into (embedded) right-angled regular polyhedra (dodecahedra or ideal octahedra) embeds geodesically in a complete finite-volume connected orientable hyperbolic 4-manifold $W$, which is also tessellated into right-angled regular polytopes (120-cells and ideal 24-cells). If $M$ is connected, then $\Vol(W) < 2^{49}\Vol(M)$. 

This applies for instance to the Borromean links complement. As a consequence, the Borromean link complement bounds geometrically a hyperbolic 4-manifold.
\end{abstract}

\maketitle

\section*{Introduction}

In this note we addess the following question.

\begin{quest} \label{main:quest}
Given a complete finite-volume hyperbolic $n$-manifold $M$, is there a connected complete finite-volume orientable hyperbolic $(n+1)$-manifold $W$ that contains $M$ as a geodesic hypersurface?
\end{quest}

If the answer is positive, we say that $M$ \emph{embeds geodesically}. We note that $W$ is assumed to be connected and orientable, but it makes perfectly sense to ask Question \ref{main:quest} for disconnected and/or non-orientable hyperbolic manifolds $M$.

Embedding geodesically a given hyperbolic manifold $M$ is not a trivial task: among the uncountably many connected orientable hyperbolic surfaces, only countably many embed geodesically, and they form a dense subset of Teichm\"uller space, as proved by Fujii and Soma \cite{FS}.

In dimension $n\geqslant 3$ we are not aware of any single $M$ which does not embed geodesically. On the other hand, only few explicit hyperbolic finite-volume $3$-manifolds $M$ are known to embed geodesically, so the question is still wide open.\footnote{Added: many examples were produced after the first draft of this preprint in \cite{KRS}.}

The main object of this paper is to provide new examples. Recall that there are precisely two right-angled hyperbolic platonic solids: the (non-compact) ideal octahedron and the (compact) right-angled dodecahedron. Inspired by \cite{tetrahedral}, we say that a complete (possibly disconnected) finite-volume hyperbolic 3-manifold is \emph{dodecahedral} or \emph{octahedral} if it may be tessellated into right-angled dodecahedra or ideal octahedra, respectively. Dodecahedral manifolds are compact, while octahedral ones have cusps. 

\begin{figure}
 \begin{center}
  \includegraphics[width = 4 cm]{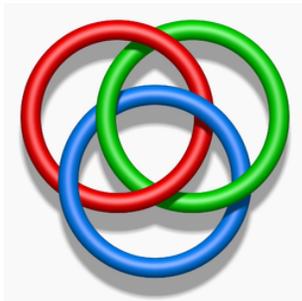}
 \end{center}
 \caption{The Borromean rings complement is hyperbolic. It is nicely octahedral and hence embeds geodesically in a finite-volume hyperbolic four-manifold.}  \label{Borromean:fig}
\end{figure}

For instance, the Whitehead and Borromean link complements (see Fig.~\ref{Borromean:fig}) are octahedral manifolds and Thurston's first famous examples of closed hyperbolic manifolds fibering over $S^1$ are dodecahedral, see \cite{Su}. 

We will need a couple of technical requirements. We say that a dodecahedral manifold is \emph{nice} if it can be decomposed into embedded dodecahedra. This is equivalent to ask that every vertex in the decomposition is adjacent to 8 distinct dodecahedra. We define an analogous but weaker requirement for the octahedral manifolds. We say that an octahedral manifold is \emph{nice} if every ideal triangle of the decomposition is adjacent to two distinct ideal octahedra (we do not require however that the four octahedra incident at every edge be distinct, so the octahedra need not to be embedded). For instance, the Whitehead link complement is not nice, but the Borromean rings complement is nice.

We prove here the following.

\begin{teo} \label{finite:teo}
Every (possibly disconnected) nice dodecahedral or nice octahedral hyperbolic three-manifold $M$ embeds geodesically in some connected complete finite-volume orientable hyperbolic four-manifold $W$.
\end{teo}

In particular, the Borromean rings complement embeds geodesically.
We note that $M$ may be disconnected: for instance, we can embed multiple copies of the Borromean rings complement disjointly in a single connected complete finite-volume hyperbolic four-manifold $W$. 

We can also estimate the volume of $W$ in terms of that of $M$. For simplicity we restrict ourselves to connected manifolds $M$. 

\begin{teo} \label{explicit:teo}
If $M$ is connected, there is a $W$ with
$$\Vol(W) \leq 13\cdot 2^{45}\cdot\Vol(M), \qquad \Vol(W) \leq 1844\cdot\Vol(M)$$
in the dodecahedral and octahedral case respectively.
\end{teo}

In our construction the hyperbolic four-manifold $W$ is tessellated into four-dimensional right-angled hyperbolic regular polytopes, the \emph{24-cell} and \emph{120-cell}, whose facets are octahedra and dodecahedra. Some of the techniques we use are borrowed from \cite{KMT}.
We prove Theorems \ref{finite:teo} and \ref{explicit:teo} in Section \ref{proof:section} and then make further comments in Section \ref{comments:section}. 

\subsection*{Acknowledgements} 
The author warmly thanks Steven Tschantz for drawing and allowing him to use the pictures in Fig.~\ref{fig5alt1:fig} and \ref{fig5balt1:fig}.
The pictures in Fig.~\ref{Borromean:fig} and Fig.~\ref{cells:fig} are taken from Wikipedia Commons: that in Fig.~\ref{Borromean:fig} lies in the Public Domain, those in Fig.~\ref{cells:fig} were produced using the software {\tt Stella} by its author \cite{Webb}. 

\section{Manifolds with right-angled corners} \label{proof:section}

We prove here Theorems \ref{finite:teo} and \ref{explicit:teo}.



\subsection{Manifolds with right-angled corners}
We now generalize both hyperbolic manifolds and right-angled polyhedra in a single notion. 

We visualize hyperbolic space via the disc model $D^n$ and define $P\subset D^n$ as the intersection of $D^n$ with the positive sector $x_1,\ldots, x_n\geqslant 0$. A \emph{hyperbolic manifold with (right-angled) corners} is a topological $n$-manifold $M$ with an atlas in $P$ and transition maps that are restrictions of isometries.

The boundary $\partial M$ is stratified into vertices, edges, $\ldots$, and facets. Distinct strata of the same dimension meet at right-angles. Examples of such $M$ are hyperbolic manifolds with geodesic boundary and right-angled polytopes.

A simple albeit crucial property is that if we glue two hyperbolic manifolds with corners $M_1$ and $M_2$ along two isometric facets, the result is a new hyperbolic manifold with corners. 

We prove here the following.

\begin{prop} \label{embeds:corner:prop}
Every octahedral or dodecahedral hyperbolic 3-manifold $M$ embeds geodesically in the interior of a connected, complete, finite-volume orientable hyperbolic four-manifold $W$ with corners.
\end{prop}

\begin{proof}
We first consider the case $M$ is orientable. The manifold $M$ is tessellated into right-angled octahedra or dodecahedra: by placing a right-angled 24-cell or 120-cell ``above'' each octahedron or dodecahedron of the tessellation, we obtain a hyperbolic four-manifold $W'$ with corners, whose boundary $\partial W'$ contains $M$ as a connected component. By doubling $W'$ along $M$ we get a hyperbolic manifold with corner $W$ containing $M$ in its interior (see Fig.~\ref{W:fig}).

\begin{figure}
 \begin{center}
  \includegraphics[width = 9 cm]{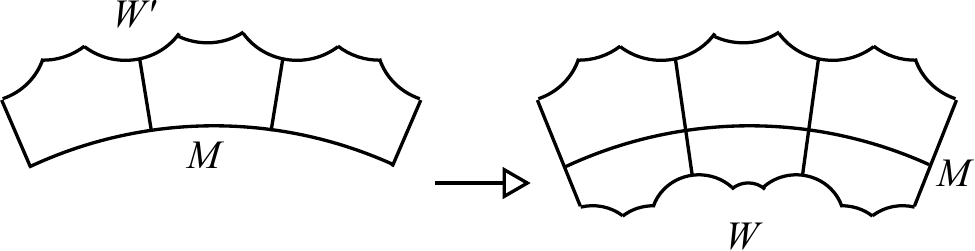}
 \end{center}
 \caption{We place a right-angled 24- or 120-cell above each octahedron or dodecahedron in $M$ to get $W'$, then we double it to get a manifold with corner $W$ containing $M$ in its interior.}
 \label{W:fig}
\end{figure}

If $M$ is connected, then $W$ also is and we are done. Otherwise, each component $M_i$ of $M$ is contained in a component $W_i$ of $W$, for $i=1,\ldots,k$. We note that every 24-cell or 120-cell has a facet opposite to that contained in $M$ which is still a facet of $W$, hence $\partial W_i$ contains at least one octahedral or dodecahedral facet $f_i$ for each $i=1,\ldots,k$.

\begin{figure}
 \begin{center}
  \includegraphics[width = 11 cm]{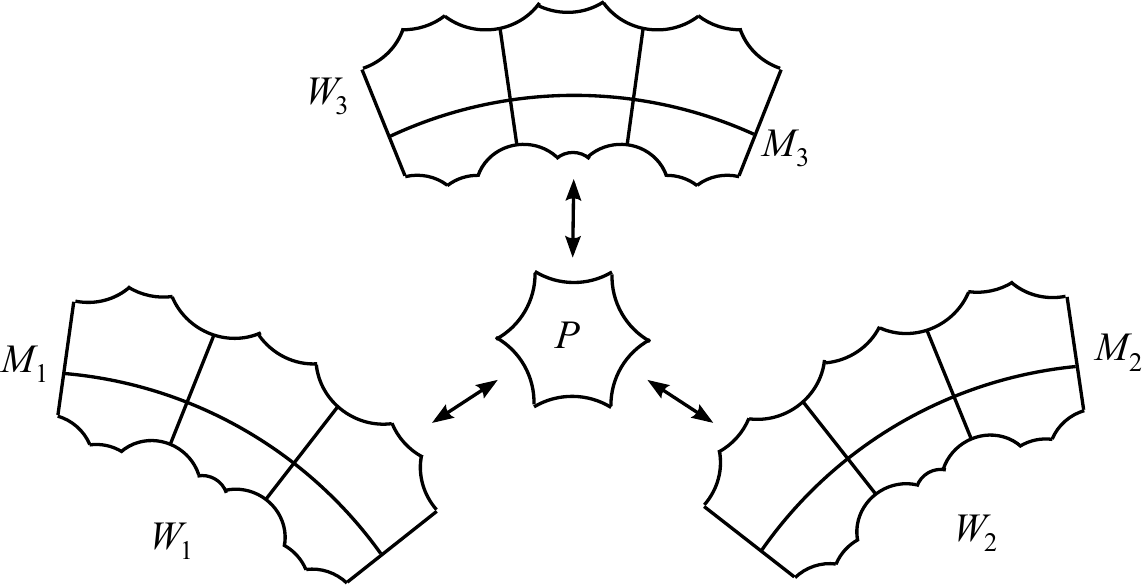}
 \end{center}
 \caption{If $M$ is not connected, we use another right-angled polytope $P$ to connect all the components $W_i$.}
 \label{W2:fig}
\end{figure}

We only need to connect all these facets $f_i$ altogether. To do this, we observe that both the 24- and the 120-cell contain (at least) three pairwise non-incident distinct facets, and by doubling the polytope multiple times along one of these we get a bigger right-angled (not regular) polytope $P$ with disjoint facets $g_1, \ldots, g_k$ that are isometric to the $f_i$. We attach $P$ to the disconnected $W$ by identifying $f_i$ to $g_i$ for every $i=1,\ldots,k$ and thus get a connected manifold with corners that we still name $W$ (see Fig.~\ref{W2:fig}).

Finally, if a component $M_i$ is non-orientable, we consider its orientable double covering $\tilde M_i$ and construct a $\tilde W_i$ containing $\tilde M_i$ as above. The orientation-reversing deck involution $\iota \colon \tilde M_i \to \tilde M_i$ extends uniquely to an orientation-\emph{preserving} involution $\iota\colon \tilde W_i \to \tilde W_i$ that exchanges the two sides of $\tilde W_i \setminus \tilde M_i$ and we set $W_i = \tilde W_i/_\iota$.
\end{proof}

\begin{oss} \label{2n:oss}
If $M$ is connected and tessellated into $n$ right-angled dodecahedra or ideal octahedra, the proof of Proposition \ref{embeds:corner:prop} produces a $W$ that is tessellated into $2n$ right-angled 120-cells or ideal 24-cells. 
\end{oss}

\subsection{Colorings}
We now show how to promote a hyperbolic manifold with corners to a manifold without corners. We first note that a manifold with corners $W$ is naturally a hyperbolic orbifold, and as such it has a finite cover that is a hyperbolic manifold by Selberg's Lemma. However, we want a manifold cover that still contains $M$, and this is not guaranteed.

We now construct some explicit manifold finite covers, following the \emph{coloring} technique used by various authors in similar contexts, see for instance Vesnin \cite{Vesnin87, Vesnin}, Davis and Januszkiewicz \cite{DJ}, Izmestiev \cite{I}, and Kolpakov, Martelli, and Tschantz  \cite{KMT}.

Let $W$ be a hyperbolic manifold with corners. Consider a finite set $\calS = \{1,\ldots,k\}$ of colors. Let a \emph{coloring} $\lambda$ of $W$ be the assignment of a color $c\in \calS$ at every facet of $W$, such that adjacent facets have different colors. In particular, the $n$ facets incident to a vertex all have distinct colors, hence if $\partial W$ contains vertices we must have $k\geqslant n$. We suppose for simplicity that all the $k$ colors in $\calS$ are used by $\lambda$ (if not, just take a smaller $k$).

\begin{prop} \label{promote:prop}
Let $W$ be a connected orientable hyperbolic manifold with corners. If $W$ has a coloring $\lambda$ with $k$ colors, there is a connected orientable hyperbolic manifold $\tilde W$ that is tessellated into $2^k$ identical copies of $W$. The induced map $\tilde W \to W$ is well-defined and is an orbifold cover of degree $2^k$.
\end{prop}
\begin{proof}
We construct $\tilde W$ as follows. We consider the $\matZ_2$-vector space $\matZ_2^k$, with canonical basis $e_1,\ldots, e_k$ and finite cardinality $2^k$. For every vector $v\in \matZ_2^k$ we define a copy $vW$ of $W$, so that we get $2^k$ disjoint identical copies overall. For a facet $F$ of $W$, we indicate by $vF$ the corresponding facet in $vW$.

For every $v\in \matZ_2^k$ and every facet $F$ of $W$, we identify $vF$ with $(v+e_c)F$ where $c$ is the color of $F$. Since adjacent facets have distinct colors, one sees easily that the result of this gluing is a genuine hyperbolic $n$-dimensional manifold $\tilde W$ without corners. If $W$ is connected and orientable, then $\tilde W$ is. The manifold $\tilde W$ is tessellated into $2^k$ copies of $W$ and we get an orbifold covering $\tilde W \to W$ of degree $2^k$.
\end{proof}

We can now refine Proposition \ref{embeds:corner:prop} by estimating the number of colors needed for $W$. The number 43 is probably not optimal, but the important point is that it does not depend on $M$. We need here the niceness hypothesis to show that there is at least one colouring.

\begin{lemma} \label{colors:lemma}
Let $M$ be a nice dodecahedral or a nice octahedral manifold. If $M$ is connected, the boundary $\partial W$ of the manifold $W$ constructed in Proposition \ref{embeds:corner:prop} can be colored with at most $43$ or $8$ colors, depending on whether $M$ is dodecahedral or octahedral.
\end{lemma}
\begin{proof}
The proof depends heavily on the combinatoric of the 120- and 24-cells, shown in Fig.~\ref{cells:fig}. We may suppose to simplify notations that $M$ is orientable: in the non-orientable case the proof is just the same. Let us first consider the case where $M$ decomposes into dodecahedra.

\begin{figure}
 \begin{center}
  \includegraphics[width = 6 cm]{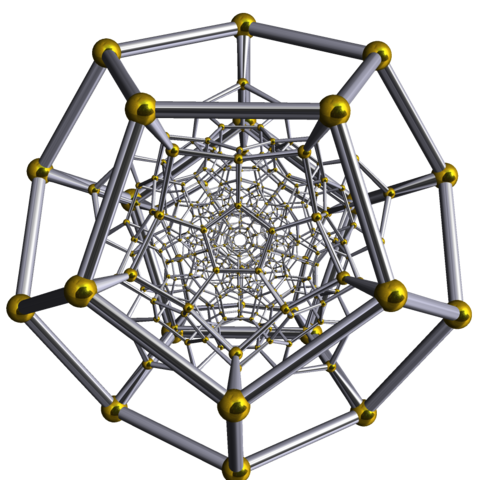}
  \includegraphics[width = 6 cm]{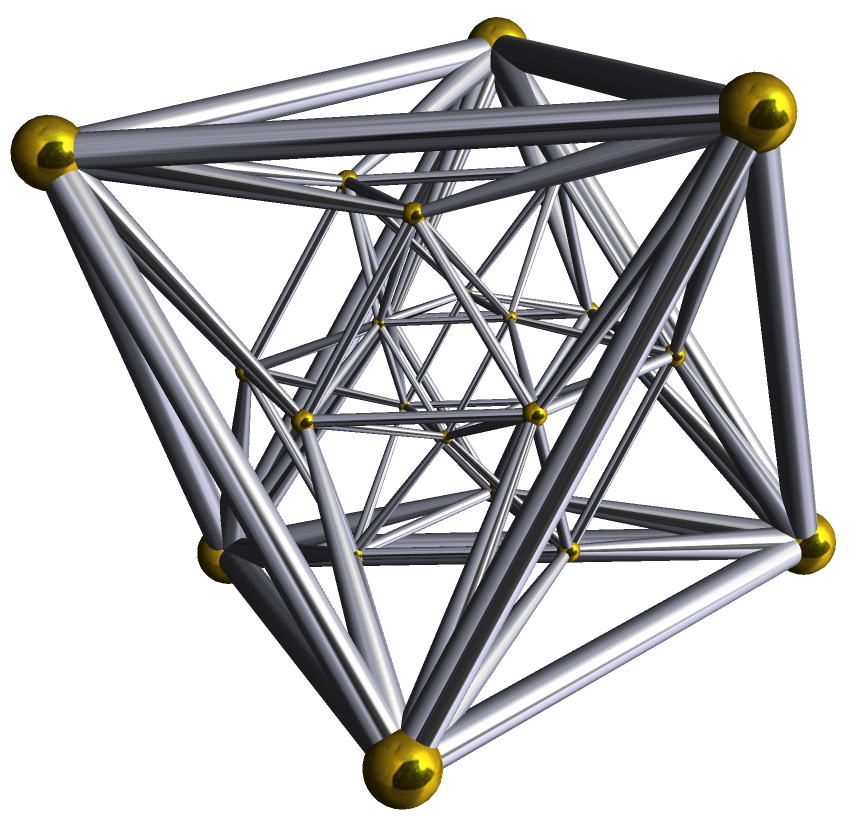}
 \end{center}
 \caption{The 120- and 24-cells. The picture shows the tessellation of $S^3$ into 120 or 24 regular polyhedra and its layers: the first layer is the unbounded face containing $\infty$ in the picture, the adjacencies between the second and third layers can be seen in the figure with little effort.}  \label{cells:fig}
\end{figure}

A right-angled 120-cell $Z$ is attached to each right-angled dodecahedron $D$ in $M$ to form $W'$, and then $W$ is the double of $W'$ along $M$. Starting from $D$ as a ``north pole'', and ending to its opposite ``south pole'', the 120 facets of $Z$ decompose into nine spherical layers consisting of 1, 12, 20, 12, 30, 12, 20, 12, and 1 dodecahedra.

The 12 facets in the second layer are incident to $D$ and are identified to the second-layer facets of the adjacent 120-cells (attached to the dodecahedra in $M$ incident to $D$). All the facets in the higher layers form the boundary $\partial W$. Note however that $\partial W$ does not consist simply of dodecahedra, because some of them are glued together with dihedral angle $\pi$ to form more complicate facets, and we now need to control this phenomenon. To understand the problem, imagine the lower dimensional situation where one attaches right-angled dodecahedra to a surface that tessellates into right-angled regular pentagons. One checks easily that in this case the faces of $\partial W$ are right-angled pentagons and octagons, the latter partitioned into four pentagons. The situation here is slightly more complicated but analogous. 

Fig. \ref{cells:fig}-(left) shows that a second-layer facet is incident to 5 third-layer facets and 1 fourth-layer one. Each third-layered facet $F$ is incident to three second-layered ones, all incident to a single vertex $v$. The facet $F$ is attached to three similar third-layered facets of other 120-cells, which are also incident to $v$, and by repeating this we get that $F$ is contained in a facet $Q$ of $\partial W$ consisting of eight dodecahedra, all sharing the same vertex $v$ which lies in the center of $Q$. The polyhedron $Q$ is shown in Fig.~\ref{fig5alt1:fig} and has 42 faces: 24 pentagons, 12 hexagons, and 6 octagons. 

\begin{figure}
 \begin{center}
  \includegraphics[width = 7 cm]{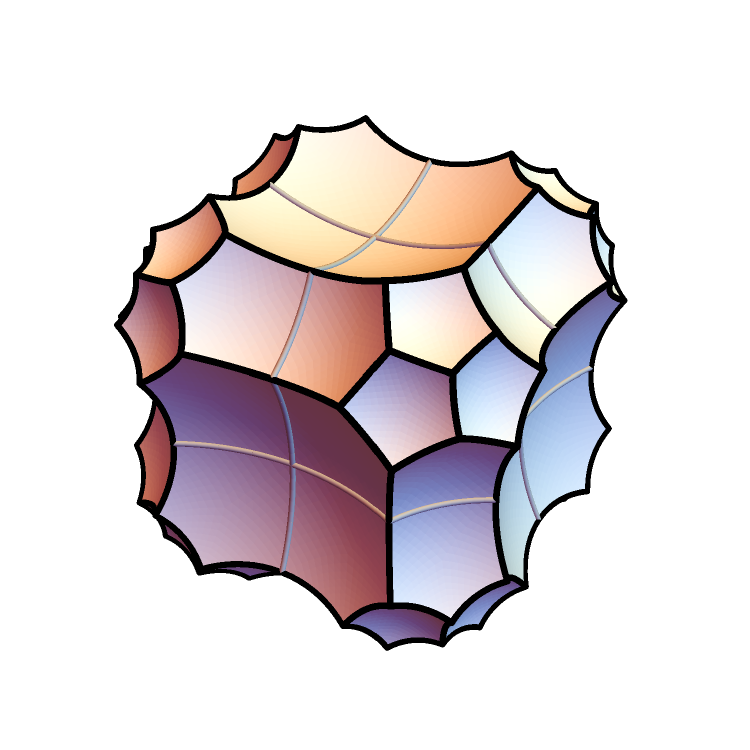}
 \end{center}
 \caption{A hyperbolic right-angled polyhedron $Q$ obtained by attaching $8$ right-angled dodecahedra to the same vertex $v$, which is the barycenter of $Q$. The polyhedron $Q$ has 42 faces: 24 pentagons, 12 hexagons, and 6 octagons.}  \label{fig5alt1:fig}
\end{figure}

\begin{figure}
 \begin{center}
  \includegraphics[width = 7 cm]{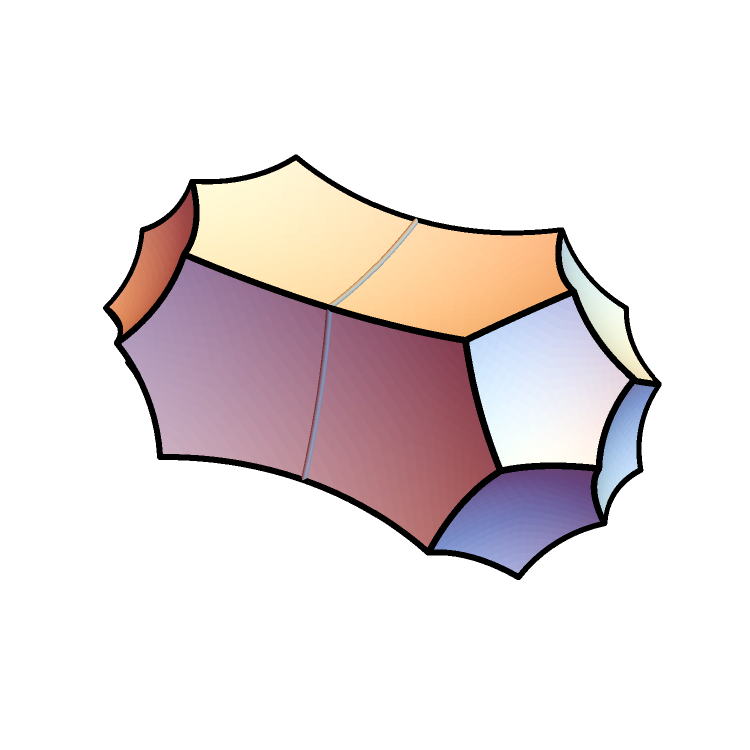}
 \end{center}
 \caption{A hyperbolic right-angled polyhedron $R$ obtained by mirroring a right-angled dodecahedron along a face.}  \label{fig5balt1:fig}
\end{figure}

Every fourth-layered facet $F$ is incident to a single second-layered one, and is hence attached to a single fourth-layered one in some adjacent 120-cell: the two form altogether the polyhedron $R$ shown in Fig.~\ref{fig5balt1:fig}, which has 17 faces: 12 pentagons and 5 hexagons. 

Summing up, the boundary $\partial W$ is tessellated into polyhedra of three types $Q$, $R$, and $D$, with 42, 17, and 12 faces respectively. The polyhedra of type $Q$ correspond to the vertices of the tessellation of $M$, those of type $R$ correspond to the pentagons of the tessellation. Since $M$ is nice, every vertex is incident to 8 distinct dodecahedra, and every pentagon is incident to two distinct dodecahedra. This implies that every polyhedron of type $Q$ is made of 8 dodecahedra that belong to 8 different 120-cells. Therefore $Q$ is embedded, and in particular it is not self-adjacent along its faces. The same holds for the polyhedra of type $R$ and $D$.

The latter technical fact implies that the adjacency graph of the facets of $\partial W$ is a graph without loops, where every vertex has valence $\leq 42$, and hence it can be colored with at most 43 colors. (Every finite graph without loops and with valence bounded by $k$ can be $(k+1)$-colored, simply by ordering the vertices and then assigning them colors in a sequence.) The proof is complete.

We now turn to the non-compact octahedral case, which is a bit different. A right-angled ideal 24-cell is attached to every ideal regular octahedron $O$ in $M$ to form $W'$. The 24-cell layers into 1, 8, 6, 8, 1 octahedra, and it has the remarkable property that it can be colored with three colors, say yellow, red, and blue, each assigned to 8 faces. The 8 octahedra in the odd layers have the same color, say yellow. Each even layer contains four red and four blue octahedra.

As above, the 8 octahedra in the second layer are identified to second-layer octahedra in the adjacent 24-cells. Fig.~\ref{cells:fig}-(right) shows that each second-layer octahedron is incident to 3 third-layer (yellow) octahedra and to 1 (non-yellow) fourth-layer one.

Every third-layer yellow octahedron is adjacent to four second-layer octahedra, and is therefore attached to four third-layer yellow octahedra in adjacent 24-cells. In contrast to the 120-cell case, these four adjacent octahedra may not intersect each other; each is adjacent to three more yellow octahedra, so that the facet $R$ in $\partial W$ containing them may consist of an arbitrarily big number of octahedra! Luckily, all these octahedra are yellow, so we color $R$ in yellow, and we are sure that there are no self-adjacencies between yellow octahedra. These big yellow facets in $\partial W$ are pairwise not adjacent, so we need only to color the rest with other colors.

The remaining fourth-layer octahedra are just paired to adjacent ones producing more facets. Each such facet $F$ is obtained by doubling a right-angled octahedron along a face, and it is easily seen to have 11 faces: 8 ideal triangles and three ideal quadrilaterals. Such facets $F$ correspond to the ideal triangles of the decomposition of $M$ into ideal octahedra. Since the decomposition is nice, the two octahedra composing $F$ belong to two distinct ideal 24-cells, and hence there are no self-adjacencies of $F$ along facets.

One such facet $F$ does not have a natural color, and it is adjacent to 5 yellow facets and 6 more uncolored facets isometric to $F$ (but distinct from $F$, by what just said). The adjacency graph of the uncolored facets in $\partial W$ has no loops and has valence 6, so it can be colored with 7 (non-yellow) colors. Therefore $\partial W$ can be colored with at most 8 colors overall.
\end{proof}

\subsection{Proof of the main result}
We can now prove Theorems \ref{finite:teo} and \ref{explicit:teo}.

\dimo{Theorem \ref{finite:teo}}
By Proposition \ref{embeds:corner:prop} the manifold $M$ embeds geodesically in the interior of a connected complete orientable finite-volume hyperbolic four-manifold $W$ with corners. By Proposition \ref{promote:prop} there is a connected complete orientable finite-volume hyperbolic manifold $\tilde W$ containing $W$ and hence $M$ geodesically. 
\finedimo

\dimo{Theorem \ref{explicit:teo}}
If $M$ is connected, then $W$ is tessellated into $2n$ right-angled 120-cells (ideal 24-cells), see Remark \ref{2n:oss}, and $\partial W$ colors with at most $43$ ($8$) colors by Lemma \ref{colors:lemma}. By Proposition \ref{promote:prop} the manifold $\tilde W$ is tessellated into $2^{43}\cdot 2n = 2^{44}n$ ($2^8\cdot 2n = 2^9n$) right-angled 120-cells (ideal 24-cells).

The volumes of the right-angled dodecahedron $D$, ideal octahedron $O$, 120-cell $Z$, and ideal 24-cell $C$, are
$$\Vol(D) = 4.3062 \ldots, \qquad \Vol(O) = 3.6638\ldots, $$
$$\Vol(Z) = \frac{34}3 \pi^2 = 111.8553\ldots, \qquad \Vol(C) = \frac 43 \pi^2 = 13.1594 \ldots $$
which gives
$$\frac{\Vol(Z)}{\Vol(D)} \leq 26, \qquad \frac{\Vol(C)}{\Vol(O)} \leq 3.6.$$
This implies the result since $3.6\cdot 2^9 \leq 1844$.
\finedimo

One can probably prove that a connected dodecaheral or octahedral three-manifold virtually embeds geodesically also using the subgroup separability property for right-angled polytopes stated in \cite[Theorem 3.1]{ALR}. The disconnected case and the bounds on the volumes however do not seem to follow easily from such arguments.

\section{Comments} \label{comments:section}
A related interesting question, already studied in the literature, is the following: 

\begin{quest}
Given a complete finite-volume orientable hyperbolic $n$-manifold $M$, is there a complete finite-volume orientable hyperbolic $(n+1)$-manifold $W$ with geodesic boundary isometric to $M$?
\end{quest}

Here both $M$ and $W$ can be disconnected. If the answer is positive, we usually say that $M$ \emph{bounds geometrically}. Theorem \ref{finite:teo} has the following corollaries. Given a manifold $M$, we denote by $2M$ the disconnected manifold that consists of two copies of $M$.

\begin{cor}
Let $M$ be a nice octahedral or a nice dodecahedral orientable manifold. The manifold $2M$ bounds geometrically.
\end{cor}
\begin{proof}
Theorem \ref{finite:teo} says that there is an orientable $W$ containing $M$ geodesically. By cutting $W$ along $M$ we get an orientable hyperbolic manifold with geodesic boundary $2M$.
\end{proof}

\begin{cor} \label{bounds:cor}
Let $M$ be a connected octahedral or dodecahedral orientable manifold. If $M$ has a fixed-point free orientation-reversing involution $\iota$, then $M$ bounds geometrically.
\end{cor}
\begin{proof}
The manifold $M$ embeds geodesically in an orientable $W$. By cutting $W$ along $M$ we get a $W'$ whose geodesic boundary consists of two copies of $M$. We can kill one boundary component by quotienting it with $\iota$, and the result is an orientable $W''$ with geodesic boundary $M$.
\end{proof}

\begin{cor}
The Borromean rings complement $M$ bounds geometrically.
\end{cor}
\begin{proof}
As one can check with SnapPy, the manifold $M$ is the orientable double cover of the non-orientable octahedral manifold $m128$, that is $N4_1^{1,1}$ in the Callahan -- Hildebrand -- Weeks cusped census \cite{CaHiWe}. Therefore $M$ has such an involution $\iota$.
\end{proof}

Hyperbolic manifolds that bound geometrically exist in all dimensions \cite{LR3} and the first three-dimensional examples are contained in \cite{RT1, RT2}. More examples were then constructed in \cite{KMT} with techniques similar to the ones used here. An explicit link complement was produced in \cite{S}.

The Borromean rings complement is tessellated into two ideal octahedra and has volume 7.32772$\ldots$ When the first version of this paper appeared, this was the smallest hyperbolic three-manifold known to bound geometrically: shortly after, it was shown by Slavich \cite{S2} that the figure-eight knot complement also bounds geometrically,  using the rectified simplex that was employed previously in \cite{KS}. The figure-eight knot complement is (together with its sibling) the smallest cusped orientable hyperbolic three-manifold and has volume $2.02988\ldots$  
To the best of our knowledge, the smallest compact one known has volume 68.8992$\ldots$ and is tessellated into 16 right-angled dodecahedra, see \cite{KMT}.

If $M$ bounds geometrically a manifold $W$, then of course it also embeds geodesically in the double of $W$. Although we suspect that the latter notion is stronger than the former, we do not know a single example of hyperbolic orientable manifold that embeds geodesically but does not bound geometrically.

A fundamental result of Long and Reid \cite{LR1} shows in fact that ``most'' closed hyperbolic 3-manifolds do \emph{not} bound geodesically: a closed hyperbolic 3-manifold $M$ that bounds geometrically must have integral $\eta$-invariant $\eta(M)\in \matZ$. Note that $\eta(\overline M) = -\eta(M)$ and hence a mirrorable 3-manifold (ie one admitting an orientation-reversing isometry) has vanishing $\eta$-invariant, coherently with Corollary \ref{bounds:cor}. We ask the following.

\begin{quest} Is there a nice dodecahedral 3-manifold with non-integral $\eta$-invariant?
\end{quest}

One such manifold $M$ would embed geodesically, but would not bound. In fact Theorem \ref{finite:teo} could suggest that many hyperbolic 3-manifolds embed geodesically, whereas as we said only few of them bound. 

Finally, we note that, from a more four-dimensional perspective, Theorem \ref{finite:teo} shows that 24-cell and 120-cell manifolds form a big set,  rich enough to contain geodesically all nice octahedral and nice dodecahedral manifolds. Some of these four-dimensional hyperbolic manifolds were constructed in \cite{RT0, KM}.

\end{document}